\documentclass[12pt]{amsart}
\usepackage{amssymb,enumerate}

\newcommand{\C}{\mathbb{C}}
\newcommand{\Z}{\mathbb{Z}}
\newcommand{\bG}{{\mathbf{G}}}
\newcommand{\mult}{^\times}
\newcommand{\form}{\langle \phantom{x},\phantom{x}\rangle}

\newcommand{\phm}{\phantom{-}}

\newcommand{\phinv}{\vphantom{^{-1}}}

\newcommand{\phnm}{\mspace{-11mu}}
\newcommand{\inv}{^{-1}}
\newcommand{\tpi}{{\widetilde{\pi}}}

\newcommand{\qq}{\mathsf{q}}
\newcommand{\mm}{\mathsf{m}}
\newcommand{\hh}{\mathsf{h}}
\newcommand{\nn}{\mathsf{n}}

\newcommand{\Yo}{Y_{\text{o}}}
\newcommand{\Yc}{Y_{\text{c}}}

\newcommand{\Ycc}{Y_{\text{cc}}}

\newcommand{\blah}{\gamma}  

\DeclareMathOperator{\Hom}{Hom}
\DeclareMathOperator{\End}{End}
\DeclareMathOperator{\card}{card}
\DeclareMathOperator{\Aut}{Aut}
\DeclareMathOperator{\ind}{ind}
\DeclareMathOperator{\tr}{tr}
\DeclareMathOperator{\Sp}{Sp}
\DeclareMathOperator{\GSp}{GSp}
\DeclareMathOperator{\GL}{GL}
\DeclareMathOperator{\SL}{SL}
\DeclareMathOperator{\SO}{SO}
\DeclareMathOperator{\GO}{GO}
\DeclareMathOperator{\OO}{O}
\DeclareMathOperator{\GU}{GU}
\DeclareMathOperator{\SU}{SU}
\DeclareMathOperator{\U}{U}

\newtheorem{theorem}{Theorem}[section]
\newtheorem{lemma}[theorem]{Lemma}
\newtheorem{prop}[theorem]{Proposition}
\newtheorem{thm}[theorem]{Theorem}
\newtheorem{conj}[theorem]{Conjecture}

\theoremstyle{definition}
\newtheorem*{remark}{Remark}
\newtheorem*{question}{Question}

\newcommand{\set}[2]{
	{\left\{\left.
	#1\vphantom{#2\bigl(\bigr)}\,\right|
	\,#2\right\}}}

\newcommand{\sett}[2]{\set{#1}{\text{#2}}}

\begin{document}
\ifx\href\undefined\else\hypersetup{linktocpage=true}\fi 
\title{On certain multiplicity one theorems}
\author{Jeffrey D. Adler}
\thanks{The first-named author was partially supported by the National
Security Agency (\#MDA904-02-1-0020).}

\address{
Department of Mathematics \\
The University of Akron \\
Akron, OH  44325-4002 \\ USA}
\email{adler@uakron.edu}

\author{Dipendra Prasad}
\address{
School of Mathematics\\
Tata Institute of Fundamental Research \\
Mumbai 400 005 \\
India}
\email{dprasad@math.tifr.res.in}

\date{30 September, 2004} 

\begin{abstract}
We prove several multiplicity one theorems in this paper.
For $k$ a local field not of characteristic two,
and $V$ a symplectic space over $k$,
any irreducible admissible
representation of the symplectic similitude group $\GSp(V)$
decomposes with multiplicity one when restricted to
the symplectic group $\Sp(V)$.
We prove the analogous result for $\GO(V)$ and $\OO(V)$,
where $V$ is an orthogonal space over $k$.
When $k$ is non-archimedean,
we prove the
uniqueness of Fourier-Jacobi
models for representations of $\GSp(4)$,
and the existence of such models for supercuspidal
representations of $\GSp(4)$.
\end{abstract}

\maketitle

\setcounter{tocdepth}{1}

\tableofcontents

\section{Introduction}
In this paper we prove several multiplicity one theorems.
Our initial
aim when writing this paper was to prove a multiplicity one theorem for
the restriction of an irreducible admissible representation of $\GSp(4)$
to $\Sp(4)$ for the $p$-adic case.
As is well known, such theorems are
easy consequences of the uniqueness of Whittaker models, when they exist.
But not every representation has a Whittaker model.
Our initial attempt was thus to
use the analogous concept of Fourier-Jacobi models (recalled below),
for which uniqueness
was proved by Baruch and Rallis \cite{br} for the case of $\Sp(4)$.
This
required us to extend their work from $\Sp(4)$ to $\GSp(4)$,
which became
a major exercise in itself, useful in its own right.

We now introduce some notation.
Let $k$ denote a local field not of characteristic two.
Let $V$ denote
a finite-dimensional vector space over $k$ with a non-degenerate
bilinear form $\form$ that is either symmetric or skew-symmetric.
Let $\U(V)$ denote the associated automorphism group:
$$
\sett{ g \in \Aut(V)}{$\langle gv_1,gv_2 \rangle  = \langle
v_1,v_2 \rangle$ for all $v_1,v_2 \in V$}.
$$
Let $\GU(V)$ denote the corresponding similitude group:
$$
\sett{ g \in \Aut(V)}{$\exists \lambda_g\in k\mult$,
$\forall v_1,v_2 \in V$,
$\langle gv_1,gv_2 \rangle = \lambda_g \langle v_1,v_2 \rangle$}.
$$
We also denote these groups by $\Sp(V)$ and $\GSp(V)$
(resp.\ $\OO(V)$ and $\GO(V)$) if $\form$ is skew-symmetric
(res.\ symmetric).
If the dimension of $V$ is $2n$, we write $\Sp(V)$ also as $\Sp(2n)$,
and $\GSp(V)$ as $\GSp(2n)$.

We now introduce the Fourier-Jacobi models.
Let $e_1$ be any nonzero vector in $V$.
Let $J$ be the stabilizer of $e_1$
in $\Sp(2n)$.
Then $J \cong \Sp(2n-2) \ltimes  H$ where $H$ is the $(2n-1)$-dimensional
Heisenberg group.
Let $Z \cong k$ be the center of $H$, and
$\psi: Z \rightarrow \C\mult$ a nontrivial character.
Let $\theta_\psi$ be the oscillator representation of $H$ with
central character $\psi$.
It is well known that $\theta_\psi$
can be extended to a representation of $\tilde{J} = \tilde{S}H$
with $\tilde{S}$ the two-fold metaplectic cover of $\Sp(2n-2)$.
Let $\sigma$ be an irreducible admissible genuine representation
(i.e., nontrivial on the kernel of the map from $\tilde{S}$ to
$\Sp(2n-2)$) of $\tilde{S}$. Then $\sigma \otimes \theta_\psi$
is an irreducible admissible representation of $\tilde{J}$
which, as both $\sigma$ and $\theta_\psi$ are genuine, is in fact
a representation of $J$.

\begin{remark}
Any irreducible admissible representation of $J$
on which $Z$ operates via $\psi$
is of this form.
\end{remark}

Baruch and Rallis prove the following theorem in \cite{br}.

\begin{thm}
Suppose $k$ is non-archimedean.
Let $\pi$ be an irreducible admissible representation
of $\Sp(4)$. Then for any irreducible admissible representation
$\mu$ of $J$ on which $Z$ operates via a nontrivial character,
$\dim \Hom_{J}[{\pi}, \mu] \leq 1$.
\end{thm}

Note that in their statement of the theorem, the field $k$
has characteristic zero.  However, the proof only requires
that the characteristic is not two.

We prove an analogous theorem for $\GSp(4)$.
Although our proof is modelled on the proof in \cite{br}, many details are
quite different; in particular, the proof for the `open cell'
(see \S\ref{sec:open}) is totally different.

\begin{thm}
\label{thm:fj}
Let $\pi$ be an irreducible admissible representation
of $\GSp(4)$. Then for any irreducible admissible representation
$\mu$ of $J$ on which $Z$ operates via a nontrivial character,
$\dim \Hom_{J}[{\pi}, \mu] \leq 1$.
\end{thm}

\begin{remark}
By Frobenius reciprocity, a $J$-invariant map from
$\pi$ to $\mu$ is equivalent to an embedding of $\pi$ into the induced
representation $\ind_J^{\GSp(V)} \mu$, called a Fourier-Jacobi
model of $\pi$.
\end{remark}

To be able to use the uniqueness theorem for the Fourier-Jacobi models
for $\GSp(4)$, we must show that they exist. This amounts to showing
that a certain representation of $J$
must have irreducible quotients when
it is known to be nonzero.
This turned out to be much more difficult than anticipated.
In fact, we prove only a special case here (\S\ref{sec:existence}),
and deduce the multiplicity one theorem for
restriction from $\GSp(4)$ to $\Sp(4)$ in general from this.
We hope that the problems encountered in this part of the paper
(and the way that we have handled them) will be of independent interest.

In a similar vein, i.e., by the method of ``models,'' we give a proof
of the multiplicity one theorem for the restricition of an
irreducible admissible representation of $\GL(n)$
to $\SL(n)$ due originally to
Tadi\'c~\cite{ta},
who proved it by an elaborate analysis using the full classification
of irreducible admissible representations of $\GL(n)$ (due to
Zelevinsky~\cite{ze}).

\begin{thm}
\label{thm:gln}
Any irreducible admissible representation of $\GL(n)$ decomposes with
multiplicity one when restricted to $\SL(n)$.
\end{thm}

However, after proving the multiplicity one theorem for $\GSp(4)$ by the
method of Fourier-Jacobi models, we realized that
a more general multiplicity one theorem for restriction from
$\GU(V)$ to $\U(V)$
is an easy consequence of a result in linear algebra (of classical groups),
combined with the usual formalism of Gelfand pairs adapted to
$p$-adic groups by Gelfand-Kazhdan \cite{gelfand-kazhdan}
and developed further by Bernstein-Zelevinsky \cite{bz}.
This lemma
in linear algebra, valid for any field of characteristic not 2,
says (in the symplectic case) that for
any $g$ in $\GSp(2n)$, $g$ and ${}^t\!g$ are conjugate
by an element of $\GSp(2n)$ of similitude $-1$.
Forms of this
lemma are available for all classical groups in \cite{mvw}.
The extension of this result of \cite{mvw} to the similitude group for 
the symplectic case was observed in \cite{p-selfdual}.
But for our purposes, its most precise form
given in a very recent paper of Vinroot \cite{vi} is what will be essential.

We prove the following theorem in this paper.

\begin{thm}
\label{thm:gspn}
Let $V$
be a finite-dimensional vector space over $k$ with a non-degenerate
symmetric or skew-symmetric form $\form$.
Then any irreducible admissible representation
$\tpi$ of $\GU(V)$ decomposes with multiplicity one when
restricted to $\U(V)$;
i.e., for any irreducible, admissible representation
$\pi$ of $\U(V)$,
$$
\dim \Hom_{\U(V)}[\tpi, \pi] \leq 1.
$$
\end{thm}

In \S\ref{sec:mult1}, we will give the rather simple proofs of
theorems~\ref{thm:gln} and~\ref{thm:gspn}.
We also state a conjecture concerning multiplicity one restriction
for more general groups.
From \S\ref{sec:fj-models} on,
we will assume that $k$ is non-archimedean,
and will work exclusively with the group $\GSp(4)$.
We prove theorem~\ref{thm:fj}
about the uniqueness of the Fourier-Jacobi
model for its representations.
Then we give another proof of
the multiplicity one theorem about restriction
from $\GSp(4)$ to $\Sp(4)$ from this point of view by proving the existence
of Fourier-Jacobi models in the supercuspidal case (\S\ref{sec:existence}),
and handling
nonsupercuspidal representations separately (\S\ref{sec:non-sc}).

\textbf{Acknowledgement:}
Theorem~\ref{thm:gspn} has long been expected,
at least in the rank-two symplectic case.
We thank Paul Sally for suggesting that we study it,
and for enthusiastic encouragement on several occasions.
We also thank E.~Moshe Baruch and David Soudry
for helpful correspondence.
We thank the Institute of Mathematical Sciences,
Singapore for the invitation to participate in their
special semester activity on representation theory in the summer of 2002,
where the two authors shared an office, and the collaboration on this
theorem was conceived. The first-named author also thanks
the Tata Institute for an invitation in the summer of 2003 to continue
our work there.

\section{Multiplicity one upon restriction}
\label{sec:mult1}

In this section, we (re)prove theorem~\ref{thm:gln}
and prove theorem~\ref{thm:gspn}.
We deal with the archimedean case first, since it is easy,
and from then on assume that $k$ is non-archimedean.

\subsection{Archimedean case}
We recall some general Clifford theory:

\begin{lemma}
If $G$ is a group with center $Z$, and $H$ is a normal subgroup with
$G/ZH$ a
finite cyclic group,
then any irreducible represntation of $G$
decomposes with multiplicity one when restricted to $H$.
(If $G$ is a real Lie group, then by a representation of $G$,
we mean
either a continuous representation in a Frechet space, or a
Harish-Chandra module.)
\end{lemma}

Suppose $k$ is archimedean.
Then the lemma implies both theorem~\ref{thm:gln}
and theorem~\ref{thm:gspn}, since the quotient in the lemma
has order $1$ or $2$.

Therefore, assume for the rest of this section that
$k$ is non-archimedean.

\subsection{Restriction from $\GL(n)$ to $\SL(n)$}
In this section only, there is no restriction
on the characteristic of $k$.

Our proof of theorem~\ref{thm:gln}
depends on the following theorem of Zelevinsky,
corollary~8.3 of \cite{ze}.

\begin{thm}
Let $\tpi$ be an irreducible admissible representation
of $\GL(n)$. Let $U_n$ be the group of upper-triangular unipotent
matrices in $\GL(n)$. Then there exists a character $\theta: U_n
\rightarrow \C\mult$ such that
$\Hom_{U_n}[\tpi, \theta] \cong \C$.
\end{thm}

\begin{proof}[Proof of Theorem~\ref{thm:gln}]
If
$\Hom_{U_n}[\tpi, \theta] \cong \C$, then
$\Hom_{U_n}[{\pi}, \theta] $ is also isomorphic to $\C$
for some irreducible admissible representation $\pi$ of
$\SL(n)$ which appears in the restriction of $\tpi$ with
multiplicity exactly one. Since the set of irreducible admissible
representations $\pi$ of $\SL(n)$ such that
$\Hom_{\SL(n)}[\tpi, \pi] \not = 0$
lies in a single $\GL(n)$-orbit
(for the inner conjugation action of $\GL(n)$ on $\SL(n)$,
and hence on representations of $\SL(n)$),
this completes the proof of the
theorem.
\end{proof}

\subsection{Restriction from $\GU(V)$ to $\U(V)$}
\label{sec:gelfand-pairs}

We will prove
theorem~\ref{thm:gspn}
by applying the method of Gelfand pairs:

\begin{thm}
\label{thm:gelfand-pairs}
Suppose $G$ is the group of $k$-points
of an algebraic $k$-group,
$H$ is the group of $k$-points of a closed $k$-subgroup,
and $G/H$ carries a $G$-invariant distribution.
Suppose that $\tau$ is an algebraic anti-involution of $G$
that preserves $H$
and takes each $H$-conjugacy class in $G$ into itself.
Then
\begin{enumerate}[(a)]
\item
Every $H$-invariant distribution on $G$ is $\tau$-invariant.
\item
For any irreducible, smooth representations $\tpi$ of $G$
and $\pi$ of $H$, with smooth duals  $\tpi^\vee$ and $\pi^\vee$
respectively,
let $m(\tpi,\pi)$ denote the dimension
of the space of $H$-invariant linear maps from $\tpi$ to $\pi$.
Then $m(\tpi,\pi) m(\tpi^\vee, \pi^\vee) \leq 1$.
\end{enumerate}
\end{thm}

\begin{proof}
Both parts of the theorem are due to Gelfand-Kazhdan as refined by
Bernstein-Zelevinsky. For part (a) we refer to theorem 6.13 of \cite{bz},
and for part (b) we refer to lemma 4.2 of \cite{p-trilinear}.
\end{proof}

Suppose that $G=\GU(V)$,
$H= \U(V)$,
and $\tpi$ is an irreducible, admissible representation of $G$.
From generalities, we know that
as a representation of $H$,
$\tpi$ decomposes into a finite direct sum
of irreducible representations:
$$
\tpi \cong \pi_1 \oplus \cdots \oplus \pi_\ell.
$$
Then
$$
\tpi^\vee  \cong \pi_1^\vee \oplus \cdots \oplus \pi_\ell^\vee.
$$
Thus, $m(\tpi,\pi)=m(\tpi^\vee,\pi^\vee)$
for any summand $\pi$ of $\tpi$.
Since $m(\tpi,\pi) \geq 1$,
theorem~\ref{thm:gelfand-pairs}
will imply $m(\tpi,\pi)=1$
as long as there
exists an anti-involution $\tau$ as in
the theorem.

Thus, we will have proved theorem~\ref{thm:gspn}
if we can find a suitable anti-involution $\tau$.
This is provided by the following lemmas,
which follow from the work of
Vinroot, cf.\ corollary~1 of \cite{vi}.

\begin{lemma}
\label{lem:vinroot-gsp}
Suppose $V$ is a symplectic space.
Fix $d \in \GSp(V)$ of similitude $-1$.
Let $\tau$ be the anti-involution
on $\GSp(V)$ defined by $\tau(g) = d {}^t\!gd^{-1}$.
Then for any $g \in \GSp(V)$, $g$ and $\tau(g)$ are conjugate
by an element of $\Sp(V)$.
\end{lemma}

\begin{lemma}
\label{lem:vinroot-go}
Suppose $V$ is an orthogonal space.
Let $\tau$ be the anti-involution
on $\GO(V)$ defined by $\tau(g) =  {}^t\!g$.
Then for any $g \in \GO(V)$, $g$ and $\tau(g)$ are conjugate
by an element of $\OO(V)$.
\end{lemma}

\subsection{A conjecture on multiplicity one restriction}
\label{sec:conj}
In this paper we have proved a multiplicity one theorem
for restriction from $\GU(V)$ to $\U(V)$ (where $\U(V)$
is symplectic or orthogonal),
as well as reproved a theorem (originally due to Tadi\'c)
about restriction from $\GL(n)$ to $\SL(n)$.
We note that the theorem about $\GU(V)$ has been proved
by a generality valid for all fields not of characteristic two,
whereas the theorem on
$\GL(n)$ is proved by both Tadi\'c and ourselves using
non-archimedean local fields
(as the general lemma from linear algebra that one
may wish to be true, i.e.,
for $h$ a fixed element of $\GL(n)$,
$A$ and $h{}^t\!A h^{-1}$
are conjugate via $\SL(n)$,
does not hold, as one can easily see).
There is, however, the
possibility that such a lemma holds for distributions on $\GL(n)$,
and therefore the multiplicity one theorem
can indeed by proved by the method of
Gelfand pairs, as developed by Gelfand and Kazhdan \cite{gelfand-kazhdan}.
This suggests the possibility that, just like the uniqueness
of Whittaker models, proved for all quasi-split groups, the following too
is true in this generality, and could be proved by analyzing invariant
distributions.

\begin{conj}
Let $\bG$ be a quasi-split reductive algebraic group over
a local field $k$.
Let $\widetilde{\bG}$ be a reductive algebraic group containing $\bG$
such that the
derived groups of $\bG$ and $\widetilde{\bG}$ are the same,
and such that
$\widetilde{\bG}/\bG$ is connected.
Then multiplicity one holds for restriction of irreducible admissible
representations of $\widetilde{\bG}(k)$ to $\bG(k)$.
\end{conj}

\begin{remark}
It is well known that multiplicity one is not true for
restriction from $D\mult$ to $SL_1(D)$, $D$ a division algebra over
a local field, so the quasi-splitness assumption seems necessary.
\end{remark}

\begin{remark}
One example for which the conjecture would be especially
useful is where $\widetilde{\bG}$ is a unitary group  $\U(n)$,
and $\bG = \SU(n)$.
Just like its close cousin $(\GL(n), \SL(n))$,
multiplicity one cannot be proved purely by methods of
linear algebra, but will require careful analysis of invariant
distributions.
\end{remark}

\section{Fourier-Jacobi models: Basic setup and notation}
\label{sec:fj-models}
Assume from now on that $k$ is non-archimedean.
Let
$$
j = \left(
\begin{smallmatrix}
 &   &  &  1 \\
 &   &  1 &  \\
 &  \phnm -1 &  &  \\
-1 &  &  &
\end{smallmatrix} \right),
$$
so that $\Sp(4)$ is the subgroup of $\GL(4)$
defined by
$$
{}^t\!g j g = j,
$$
and $\GSp(4)$ is the subgroup of $\GL(4)$ defined by
$$
{}^t\!gjg = \lambda(g)j \quad\text{for some $\lambda(g)\in k\mult$}.
$$
Let $C$ denote the center of $\GSp(4)$.
For $\lambda \in k\mult$, $v\in k^2$, $A\in\GL_2(k)$, $z\in k$,
and $B\in M_2(k)$ with $B_{11}=B_{22}$, let
\begin{align*}
\mm(\lambda,A)
&= \begin{pmatrix}
\lambda & & & \\
& A_{11} & A_{12} & \\
& A_{21} & A_{22} & \\
&&& \lambda\inv \det A
\end{pmatrix}
\\
\hh(v,z)
&= \begin{pmatrix}
1 & v_1 & v_2 & z \\
& 1 & & v_2 \\
&& 1 & \phnm -v_1 \\
&&& 1
\end{pmatrix}
\\
\qq(A,B,\lambda)
&=
\begin{pmatrix}
A & \\ & \lambda A^*
\end{pmatrix}
\begin{pmatrix}
I & B \\ & I
\end{pmatrix},
\\
& \qquad \text{where
$A^* = \omega {}^t\!A\inv \omega\inv$,
$\omega = \left(\begin{smallmatrix}0& 1\\1&0\end{smallmatrix}\right)$}
\\
\nn(B)
&= \qq(I,B,1).
\end{align*}
Let $M$, $H$, and $N$ denote the images of $\mm$, $\hh$,
and $\nn$, respectively.
Then $P=MH$ is the Klingen parabolic subgroup
of $\GSp(4)$, and its unipotent radical $H$ is the Heisenberg group.
The image of $\qq$ is the Siegel parabolic subgroup $Q$, whose
unipotent radical is $N$.
Let
\begin{align*}
M' &= \set{\mm(1,m)}{m\in\SL(2)} \\
J &= M'H  \quad\text{(the \emph{Fourier-Jacobi} group)}\\
L &= \set{
	\left(
	\begin{smallmatrix}
	1&&&\\&1&&\\&&\lambda&\\&&&\lambda
	\end{smallmatrix}
	\right)
	}{\lambda\in k\mult} \\
Z &= \set{\hh(0,z)}{z\in k}.
\end{align*}
Note that $Z$ is the center of both $H$ and $J$.
Note also that $M=CLM'$ (in any order), and thus
$P=JCL$.

Let $\tau$ be the involution on $\GSp(4) \times J$ defined by
$$
\tau(g,h) = (d^{-1}j^{-1}{}^t\!gjd, d^{-1}h^{-1}d),
$$
where
$$
 d = \left(
\begin{smallmatrix}
-1 &   &  &   \\
 &   &  1 &  \\
 &  1 &  &  \\
 &  &  & 1
\end{smallmatrix} \right).
$$
Clearly, $d$ normalizes the subgroup $J$, and the involution
$\tau$ when restricted to the center $Z$ of the Heisenberg
group $H$ is trivial.
We will abuse notation to denote the restriction of $\tau$ to any
$\tau$-invariant subgroup of $\GSp(4)$ also by $\tau$.
We note that
$$
\tau(\mm(\lambda, A)) = \mm(\lambda^{-1}\det A, d_1^{-1}{}^t\!Ad_1),
$$
where $d_1 = \left(
\begin{smallmatrix}
 1 &  \phm 0 \\
 0 & -1
\end{smallmatrix} \right)$.

Let $\Delta J$ denote the image of $J$ under the diagonal
embedding $J\rightarrow \GSp(4) \times J$.
By the method of Gelfand pairs, as developed by
Gelfand and Kazhdan \cite{gelfand-kazhdan}
(and applied, for example, in \cite{br} and in \S\ref{sec:gelfand-pairs}),
to prove theorem~\ref{thm:fj},
it suffices to show that any distribution on
$\GSp(4) \times J$ which is bi-invariant under $\Delta J \subset
\GSp(4)
\times J$, and $\psi$-quasi-invariant under translations of the second
variable by $Z$, is fixed by the involution $\tau$.

Just as in theorem~2.6 of \cite{br},
this is equivalent to proving the following theorem.

\begin{thm}
\label{thm:tau-inv}
Let $T$ be a distribution on $\GSp(4)$ which is invariant
under inner-conjugation by $J$, and $L_z\cdot T = \psi(z) T$
for all $z \in Z$ (where $L_z$ is left translation by $z$).
Then $T$ is fixed by $\tau$.
\end{thm}

This is clearly equivalent to the following:

\begin{thm}
\label{thm:T-zero}
Let $T$ be as in theorem~\ref{thm:tau-inv},
and suppose in addition that $T$ is $\tau$-skew-invariant.
Then $T$ is identically zero.
\end{thm}

\section{General strategy for proving uniqueness}
We outline the general strategy of our proof of the
theorem~\ref{thm:T-zero}.
Implicitly, it involves decomposing $\GSp(4)= P\cup Pw'P \cup Pw''P$
into a disjoint union of $J$-invariant, $\tau$-invariant
subsets
$X_0, X_1, \ldots, X_m=P$,
such that $Y_i= \cup_{j \geq i} X_j$ is a closed subset
of $\GSp(4)$, and $X_i$ is open in $Y_i$.
We begin by showing that $T$
vanishes on the open subset $X_0$, and thus restricts
to its complement.  Continuing in this way, we will
show in turn that $T$ vanishes on the complement
of $Y_i$ for all $i$
(the final case, i.e., vanishing of $T$ on $P$
is the subject of \S\ref{sec:closed}), completing the proof that $T=0$.
We emphasize that the method used in the $i$th step will vary
with $i$.
In many cases, we will show that every $J$-orbit in $X_i$
is $\tau$-stable, and use the following lemma of Bernstein,
cf.\ lemma~2.7 of \cite{br}.

\begin{lemma}
\label{lem:bernstein}
Let $X$ be the set of $k$-points of a $k$-variety
on which a group $J$ acts, as well as an
automorphism $\tau$
of order two normalizing the action of $J$, i.e., in the automorphism
group of $X$, $\tau J \tau^{-1} = J$.
If every $J$-orbit in $X$ is stable under $\tau$, then
every $J$-invariant distribution on $X$ is $\tau$-invariant.
\end{lemma}

In some cases we will show that every  $J$-orbit in $X_i$ is
stable under left multiplication by $Z$, and appeal
to the following lemma, cf.\ lemma~2.8 of \cite{br}.

\begin{lemma}
\label{lem:Z-stable}
Let $X$ be a $J$-stable subvariety of $\GSp(4)$ which is stable under $Z$
(where $J$ acts by conjugation and $Z$ by right
translation).
If every $J$ orbit in $X$ is stable under $Z$, then a distribution
on $X$ which is $J$-invariant, and on which $Z$ operates by $\psi$,
is trivial.
\end{lemma}

These properties of the orbits imply
that all distributions on $X_i$ (not just
those on the closure of $X_i$) with our invariance properties
must vanish on $X_i$, a stronger result than we need.
In a few cases, we will have to use more delicate means
to show that $T$ vanishes on $X_i$, but these cases can be reduced
to \cite{br}, which is what we do in this paper.

\section{Using the result of Baruch and Rallis on $\Sp(4)$}
\label{sec:br}

Let ${\mathcal G} = k\mult  \Sp(4)$.
Clearly ${\mathcal G}$ is an open subgroup of
$\GSp(4)$, and therefore any distribution on $\GSp(4)$
can be restricted to it. Let $T$ be a
distribution on $\GSp(4)$
with invariance properties under $J$ and $Z$ as in the statement
of theorem~\ref{thm:T-zero}, and which
transforms under $k\mult$ by a given character
(the central character), and is $\tau$ skew-invariant.
The restriction to ${\mathcal G}$ of
such a distribution
is equivalent (by a form of Frobenius reciprocity)
to a distribution on $\Sp(4)$
with invariance under $J$ and $Z$ and which further
is $\tau$ skew-invariant.
The $\Sp(4)$ theorem of Baruch and Rallis implies that this
distribution on $\Sp(4)$ is zero.
Hence our distribution $T$ is zero on this subgroup ${\mathcal G}$.
In the next sections, we analyze the possible support for the
distribution $T$.

\section{Open cell}
\label{sec:open}

This section is devoted to proving the following result:

\begin{lemma}
\label{lem:open-cell}
Every distribution on $Pw''P$ satisfying the invariance properties
of theorem~\ref{thm:tau-inv} is $\tau$-invariant.
\end{lemma}

\subsection{Transferring the problem from $Pw''P$ to smaller spaces.}
Let $X = Pw''P$.  Clearly,
$X = HMw''H$, with $H$ acting on $X$ by conjugation.
Thus $H$-invariant distributions on
$$
X = HMw''H \cong HMH
$$
can be identified with distributions on
$P=MH$ under the map
$$
h_1mw''h_2 \in HMw''H \mapsto mh_2h_1 \in MH.
$$
Since $w''$ commutes with $M' \cong \SL(2)$,
for elements $m_1 \in \SL(2)$,
$$
m_1(h_1mw''h_2) m_1^{-1} = (m_1h_1m_1^{-1})
(m_1mm_1^{-1})w''(m_1h_2m_1^{-1}).
$$
Therefore under the identification of $H$-invariant distributions
on $X = HMw''H$ with distributions on $MH$, the $J$-invariant
distributions correspond to distributions on $M \times H$ on which
$\SL(2)$ operates in the natural way by the inner-conjugation action.

It can be checked that $\tau(w'') = w''$. Therefore,
for $g = h_1mw''h_2$, $\tau(g) = \tau(h_2) w'' \tau(m) \tau(h_1)
= \tau(h_2)w''\tau(m) w''^{-1} w'' \tau(h_1)$.
Therefore under the identification of distributions on
$X = HMw''H$ with distributions on $MH$ through the
map $(h_1mw''h_2) \mapsto (m,h_2h_1)$,
the involution $\tau$ on $X$ corresponds to the involution
$(m,h) \mapsto (w''\tau(m)w''^{-1},\tau(h))$.
Thus we are reduced to proving that
$\SL(2)$-invariant distributions on $M \times H$ are invariant under
this latter involution.

Actually, we are looking at distributions on $X = HMw''H$
on which $Z$ acts on the left via $\psi$.
Clearly,
distributions on  $X$ which are $H$-invariant and $(Z,\psi)$-invariant
correspond to distributions on $MH$ which
are $(Z,\psi)$-invariant.
These correspond to distributions
on $MH/Z \cong M \times k^2$.

Since $w''\tau(m) w''^{-1} = \mm(\lambda,d_1^{-1}{}^t\!Ad_1)$,
for $m = \mm(\lambda,A)$,
we are finally reduced to proving
the following result:

\begin{lemma}
\label{lem:gl2k2}
An $\SL(2)$-invariant distribution on $\GL(2) \times k^2$
is invariant under
$\tau':(g,v) \mapsto (d_1^{-1}{}^t\!gd_1, d_2v)$,
where
$d_1 = \left(
\begin{smallmatrix}
1 & \phm 0 \\
0 & -1
\end{smallmatrix} \right)$
and
$d_2 = \left(
\begin{smallmatrix}
0 & 1 \\
1 & 0
\end{smallmatrix} \right)$.
\end{lemma}

Before we proceed further, we note the following lemma.

\begin{lemma}
\label{lem:sl2-conj}
For any $g,w \in \GL(2)$
with $\det w = -1$, the matrices
$g$ and $w{}^t\!gw^{-1}$
are conjugate by an element of $\SL(2)$.
\end{lemma}

This is a special case of lemma~\ref{lem:vinroot-gsp}
(and is not difficult to prove directly).

\begin{remark}
It follows from this lemma that representations of
$\GL(2)$ restrict to $\SL(2)$ without multiplicity,
something that is already
clear from Whittaker model considerations.
Unfortunately, there is no analogue
of lemma~\ref{lem:sl2-conj} for higher $n$, and therefore
there is no Gelfand pairs proof of Tadi\'c's theorem.
\end{remark}

\subsection{On a certain quadratic form}

Let $k^2$ be the 2-dimensional vector space over $k$ with the
standard symplectic structure $\form$. Associated to any $g \in \GL_2(k)$,
we have a quadratic form $Q_g$ on $k^2$ defined by
$$
Q_g(v) = \langle gv,v \rangle.
$$
It can be seen that $Q_g$ is a non-degenerate quadratic form
if and only if the eigenvalues of $g$ (in the algebraic closure $\bar{k}$
of $k$) are distinct. However, we will not have any occasion to use
this fact.

\begin{lemma}
\label{lem:soQ}
For an element $g \in \GL(2)$, let $Z(g)$
denote its centralizer in $\GL(2)$. Then we have $\SO(Q_g) = Z(g) \cap
\SL(2)$.
\end{lemma}

\begin{proof}
Clearly,
\begin{align*}
t \in \SO(Q_g)
& \iff  Q_g(tv) = Q_g(v)
&& \text{for all $v \in k^2$, and $\det t = 1$},\\
& \iff  \langle g tv,tv \rangle = \langle gv,v \rangle
&& \text{for all $v \in k^2$, and $\det t = 1$},\\
&\iff  \langle t^{-1}gtv,v \rangle = \langle gv,v \rangle
&& \text{for all $v \in k^2$, and $\det t = 1$},\\
& \iff  \langle [g - t^{-1}gt]v,v \rangle = 0
&& \text{for all $v \in k^2$, and $\det t = 1$}.
\end{align*}
Observe that for $v \not = 0$, $\langle w,v \rangle = 0$ if and only if
$w = \lambda v$ for some $\lambda \in k$. Therefore,
an element $t \in \SL(2)$ belongs to $\SO(Q_g)$ if and only if
for any $v \in k^2$,
$[g - t^{-1}gt]v = \lambda_vv$ for some $\lambda_v \in k$.

It is well known that if every vector of a vector space is an eigenvector
for a given linear operator, then the linear operator must be a multiple
of the identity. Therefore,
$$
g - t^{-1}gt  = \lambda I, \quad\text{for some $\lambda \in k$}.
$$
Taking the trace, we find that $\lambda$ must be zero, i.e., $t \in
Z(g)$. The lemma follows.
\end{proof}

\begin{lemma}
\label{lem:soQ-conj}
For any quadratic form $q$ on $k^2$ and vectors $v_1,v_2$
with $q(v_1) = q(v_2) \not = 0$, there exists an element $
g \in \SO(q)$ with $gv_1 =v_2$.
\end{lemma}

\begin{proof}
This is the usual Witt's theorem, except for the conclusion
that $g$ can be chosen to have determinant 1. Given $v_1,v_2$ with
$q(v_1) = q(v_2) \not = 0$, there exist $w_1 \perp v_1$ and
$w_2 \perp v_2$.
Since the discriminant of $q$ is equal to
$q(v_1)q(w_1)$, as well to $q(v_2)q(w_2)$,
we may assume that that $q(w_1) = q(w_2)$.
Clearly, the transformations that take $v_1$ to $v_2$ and $w_1$ to
$\pm w_2$ are in $\OO(q)$, and one of them has determinant 1.
\end{proof}

\subsection{Proof of Lemma~\ref{lem:gl2k2}}
Let
$$
\begin{matrix}
\pi:  &  \GL(2) \times k^2  & \rightarrow  & k \times k \times k \\
 &   &   &  \\
 & (g,v) & \mapsto & (\tr(g), \det(g), \langle gv,v \rangle),
\end{matrix}
$$
where $\form$ is the standard symplectic form on $k^2$ with
$\langle e_1,e_2 \rangle =1= -\langle e_2,e_1 \rangle$, and $\langle
e_1,e_1 \rangle = \langle e_2,e_2\rangle =0$.

It is easy to see that this mapping is $\SL(2)$-invariant and is also
$\tau'$-invariant. (For $\tau'$ invariance, we note that
$\langle gv,w\rangle = \langle v,d''{}^t\!gd''^{-1}w \rangle$, for
$d'' = \left(
\begin{smallmatrix}
\phm 0 &   1 \\
-1 & 0
\end{smallmatrix} \right)$,
and that $d''= d_1d_2$.)
We will prove the proposition by showing
that any $\SL(2)$-invariant distribution supported
on a fiber of $\pi$ is $\tau'$-invariant.
This is sufficient by the Bernstein Localization theorem
(lemma~\ref{lem:bernstein}).
We will achieve this
by dividing the possible fibers into three cases.
But first we introduce the following notation.

If two elements $(g_1,v_1)$ and $(g_2,v_2)$ are in the same
$\SL(2)$-orbit, i.e., there exists $s \in \SL(2)$ such that
$(g_2,v_2) = (sg_1s^{-1}, sv_1)$, we write
$(g_1,v_1) \sim_{\SL(2)} (g_2,v_2)$.

{\bf Case 1:}
Consider a fiber of $\pi$ lying over $(x,y,z)$, where $z\neq 0$.
We will show any $\SL(2)$-orbit
in such a fiber is $\tau'$-invariant.
That is, we will prove that
$$
(g,v) \sim_{\SL(2)} (d_1^{-1}{}^t\!gd_1, d_2v),
$$
for any $(g,v)$ such that
$\langle gv, v \rangle \neq 0$
where
$d_1 = \left(
\begin{smallmatrix}
1 & \phm 0 \\
0 & -1
\end{smallmatrix} \right)$
and
$d_2 = \left(
\begin{smallmatrix}
0 & 1 \\
1 & 0
\end{smallmatrix} \right)$.

By lemma~\ref{lem:sl2-conj}, there exists $s \in \SL(2)$ such that
$ d_1^{-1}{}^t\!gd_1 = sgs^{-1}$.
Further by combining
lemmas~\ref{lem:soQ} and~\ref{lem:soQ-conj},
there exists $ t \in Z(g) \cap \SL(2)$ such that $s^{-1}d_2 v = tv$.
Therefore,
$$
\begin{array}{rll}
(d_1^{-1}{}^t\!gd_1, d_2v) & = & (sgs^{-1},d_2v) \\
& \sim_{\SL(2)}&  (g,s^{-1}d_2v) \\
& = & (g, tv) \\
& \sim_{\SL(2)}& (g,v).
\end{array}
$$

{\bf Case 2:}
We next look at the fiber over an element $(x,y,z)$ with $x^2 \not = 4y$
and $z = 0$. Since $z =0$, for an element $(g,v)$ in the fiber,
$\langle gv,v \rangle = 0$, and therefore $v$ is an eigenvector of $g$.
Since $x^2 \not = 4y$, eigenvalues of $g$ (in $\bar{k}$) are distinct,
hence $g$ is diagonalizable over $k$ with distinct eigenvalues. Call
such a fiber $F_{(x,y,z)}$. We have a map
$$
\begin{matrix}
\mu:  &  F_{(x,y,z)}  & \rightarrow  &  \GL(2) \\
 &   &   &  \\
 & (g,v) & \mapsto & g.
\end{matrix}
$$

Since $g$ is diagonalizable over $k$, any conjugate of $g$
by $\GL(2)$ is in fact conjugate by $\SL(2)$. Therefore the
image of $F_{(x,y,z)}$ in $\GL(2)$ is a homogeneous space
for the $\SL(2)$ action, and can be taken to be $\SL(2)/T$
where $T$ is the subgroup of $\SL(2)$ consisting of those elements
that commute with $g$. We assume that
$g = \left(
\begin{smallmatrix}
 \alpha &   0 \\
0 & \beta
\end{smallmatrix} \right)$.
Therefore $\mu^{-1}(g) \cong \set{ (v_1,v_2) \in k^2}{ v_1  v_2 = 0}$.
Since $\tau'v = d_2 v$, $\tau'$ maps $(v_1,v_2)$ to $(v_2,v_1)$.

From a form of Frobenius reciprocity, as given for instance in~\cite{be},
cf.\ lemma on page 60, the $\SL(2)$-invariant
distributions on $F_{(x,y,z)}$
are in natural correspondence with the $T$-invariant distributions on
$\mu^{-1}(g) \cong \set{(v_2,v_1) \in k^2}{v_1 \cdot v_2 = 0}$,
where $T$ is the diagonal subgroup of $\SL(2)$ which acts on
$\mu^{-1}(g)$ by $t\cdot(v_1,v_2) = (tv_1,t^{-1}v_2)$. The following lemma
therefore suffices to prove that any
$\SL(2)$-invariant distribution on $F_{(x,y,z)}$ is $\tau'$-invariant.
This simple and basic lemma has appeared in many people's works
on invariant distributions; we refer to lemma~4.6 of \cite{p-trilinear}.

\begin{lemma}
Let $X = \set{(v_1,v_2)}{ v_1 v_2 = 0} \subset k^2$.
Let $k\mult$ operate on $X$ by $t\cdot(v_1,v_2) = (tv_1,t^{-1}v_2)$.
Then any distribution
on $X$ which is invariant under $k\mult$ is invariant under the involution
$(v_1,v_2) \mapsto (v_2,v_1)$.
\end{lemma}

{\bf Case 3:}
We finally look at the fiber over an element $(x,y,z)$ with $x^2 = 4y$
and $z =0$. We assume without loss of generality that $(x,y) = (2,1)$,
so that we are dealing with unipotent matrices.
The fiber is thus
$$
F_{(2,1,0)} = \sett{ (g,v)}{$g$ is unipotent and $gv = v$}.
$$
In this case, we will again prove that
$$
(g,v) \sim_{\SL(2)} (d_1^{-1}{}^t\!gd_1, d_2v),
$$
for any $(g,v)$ in such a fiber, and therefore
that any $\SL(2)$-invariant distribution supported
on such a fiber is  $\tau'$-invariant.

We will
find it more convenient to check invariance under the involution
$\tau'':(g,v) \mapsto (d_2^{-1}{}^t\!gd_2, d_1v)$, which differs
from $\tau'$ by an element of $\SL(2)$.

By (the proof of) lemma~\ref{lem:sl2-conj}, we can assume that
$d_2^{-1}{}^t\!gd_2 = sgs^{-1}$ with $s = \left(
\begin{smallmatrix}
1 &   n \\
0 & 1
\end{smallmatrix} \right)$.
Therefore,
$$
(d_2^{-1}{}^t\!gd_2, d_1v)
= (sgs^{-1}, d_1v) \sim_{\SL(2)} (g,s^{-1}d_1v).
$$
We will be done if $s^{-1}d_1v = v$, or $sv = d_1v$.
Since we have $z=0$, $s^{-1}d_1v$ is in any case
an eigenvector of $g$. We will assume that $g$ is a unipotent matrix
which is not identity, as the other case is trivial. Therefore,
$g$ has a unique eigenvector up to scaling. Therefore,
$$
\lambda s^{-1}d_1v = v
$$
for some $\lambda \in k\mult$.
It suffices to prove that $\lambda$ can be taken to be 1. We write
out the equation,  $\lambda s^{-1}d_1v =  v$, or
 $sv = \lambda d_1v $, assuming that $v$ is the column vector
$(v_1,v_2)$, explicitly:
$$
 \left(
\begin{matrix}
1 &   n \\
0 & 1
\end{matrix} \right)
\left(  \begin{matrix}
v_1 \\
v_2
\end{matrix} \right) =\lambda  \left(
\begin{matrix}
1 &  \phm 0 \\
0 & -1
\end{matrix} \right)
\left(  \begin{matrix}
v_1  \\
v_2
\end{matrix} \right).
$$
Equivalently,
$$
\left(  \begin{matrix}
v_1+nv_2 \\
v_2
\end{matrix} \right) =  \lambda \left(
\begin{matrix}
\phm v_1\\
-v_2
\end{matrix} \right).
$$
Therefore if $v_2 \not = 0$, $\lambda = \pm 1$. By changing $s$ to $-s$,
we then can assume that $\lambda =1$, and we are done. If $v_2 =0$, then
again $\lambda =1$.

\section{Middle cell}
\label{sec:middle}

We will prove the following:

\begin{lemma}
\label{lem:middle-cell}
Every distribution on $P \cup Pw'P$ satisfying the invariance properties
of theorem~\ref{thm:T-zero} vanishes on $Pw'P$.
\end{lemma}

As for the open cell, we need to
examine the various $J$-orbits in $Pw'P$.
Recall that we are using $\form$ to denote the symplectic
form given by the skew-symmetric matrix
$$
j = \left(
\begin{smallmatrix}
 &   &  &  1 \\
 &   &  1 &  \\
 &  \phnm -1 &  &  \\
-1 &  &  &
\end{smallmatrix} \right).
$$
Let $\{e_1,e_2,e_3,e_4\}$ denote the standard basis of $k^4$.
With this notation, since an element $p$ of $P$ has the property that
$pe_1$ is a multiple of $e_1$, and $pe_4$ is a multiple of $e_4$,
it follows that for
$g = p_1w'p_2 \in Pw'P$,
$$
g_{41}= -\langle ge_1,e_1 \rangle
= -\langle p_1w'p_2 e_1, e_1 \rangle
= -\langle w' \lambda_2e_1,\lambda_1e_1 \rangle
= -\langle \lambda_2e_2,\lambda_1e_1 \rangle
=0,
$$
where $\lambda_1$ and $\lambda_2$ are scalars.
Clearly, $g_{41}$ is zero for elements of $P$ too.
On the other hand, it can be easily checked that $\langle ge_1,e_1
\rangle \not = 0$ for $g = p_1w''p_2 \in Pw''P$. Thus,
$P \cup Pw'P$ consists exactly of those elements $g$ of $\GSp(4)$
with $g_{41}=0$.

Next note that the function $\blah: G \rightarrow k$ defined by
$\blah(g) = \langle g^2 e_1,e_1 \rangle$ is invariant
under $J$, i.e., $\blah(tgt^{-1}) = \blah(g)$ for all $t \in J$.
For $g \in Pw'P$, since $g_{41} = 0$,
$$
\blah(g) = (g^2)_{41}=  g_{42}g_{21} + g_{43}g_{31}.
$$
It can be easily checked that
$\blah$ is invariant under the action of $\tau$.

Since $w'$ normalizes $L$, we have
$$
Pw'P = JCLw'LCJ = JCLw'J,
$$
every element of which is $J$-conjugate to an element
of $JCLw'$.
So, modulo $C$, every element of $Pw'P$ is $J$-conjugate
to an element of the form
$hm'w'$, where
$h\in H$ and $m'\in M'L$.
Write
$$
g= c \: \hh((a',b'),z') \: \mm(1,m) \: w',
$$
where $c \in k\mult$ and $m\in \GL(2)$.
Let $\lambda = \det m$.
Then $\blah(g) = c^2 m_{21}\lambda$.

\begin{lemma}
\label{lem:middle-open}
Suppose $\blah(g) \neq 0$.
Then the $J$-orbit of $g$ is $\tau$-invariant.
\end{lemma}

\begin{proof}
Without loss of generality, assume $c=1$.
Let $\blah_0 = \blah(g) = m_{21}$.
Let $r= -m_{11} / \blah_0$
and $q=  m_{22} / \blah_0$.
Then $q$ and $r$ are the unique values
so that
$$
g' :=
\nn\!\begin{pmatrix}0&q\\r&0\end{pmatrix}
\:
g
\:
\nn\!\begin{pmatrix}0&q\\r&0\end{pmatrix}\inv
=
\hh((a,b''),z'')
\:
\mm(1,\begin{pmatrix} 0 & -\lambda\blah_0\inv \\ \blah_0&0 \end{pmatrix})
\:
w'
$$
for some
$a$, $b''$, and $z''$.
Let $s= -\lambda a / \blah_0$.
Then $s$ is the unique value so that
$$
g'' :=
\nn\!\begin{pmatrix}s&0\\0&s\end{pmatrix}
\:
g'
\:
\nn\!\begin{pmatrix}s&0\\0&s\end{pmatrix}\inv
=
\hh((0,b),z)
\:
\mm(1,\begin{pmatrix} 0 & -\lambda\blah_0\inv \\ \blah_0&0 \end{pmatrix})
\:
w'
$$
for some $b$ and $z$.
In other words, $g''$ is the unique element in the $J$-conjugacy
class of $g$ having this form.

Since $\tau(g'') \in Pw'P$ and $\blah(\tau(g'')) = \blah_0$,
we see that, from the calculations above, $\tau(g'')$
(like $g$) has a unique $J$-conjugate of the form
$$
\hh((0,b'),z')
\:
\mm(1,\begin{pmatrix} 0 & -\lambda\blah_0\inv \\ \blah_0&0 \end{pmatrix})
\:
w'
$$
for some $b'$ and $z'$.
The characteristic polynomial of this element (and thus of $\tau(g'')$)
is
$$
X^4 - b' \blah_0 X^3 + z' \blah_0 X^2 - b' \blah_0 X + \lambda.
$$
One can similarly compute the characteristic polynomial of $g''$
(and thus of $g$).
But since $g$ and $\tau(g'')$ must have the same characteristic
polynomial, we must have that $b=b'$ and $z=z'$.
\end{proof}

From now on, assume that $\blah(g) = 0$.
Then $m_{21}=0$, so
we may write
$$
g = \qq(A, \left(\begin{smallmatrix}r&s\\t&r\end{smallmatrix}\right),\lambda)
\in Q \smallsetminus P
$$
for some $r,s,t\in k$ and $\lambda \in k\mult$.
Since all considerations in the rest of the section depend only on
$g$ up to scalars, we assume that
$A = \left(\begin{smallmatrix}au&1\\u&0\end{smallmatrix}\right)\in\GL_2(k)$.
Let
$\beta_1$ and $\beta_2$ denote the (generalized) eigenvalues of $A$.

\begin{lemma}
\label{lem:eigens}
Suppose that for all $i,j\in\{1,2\}$,
we have $\beta_i\beta_j \neq \lambda$.
Then the $J$-orbit of $g$ contains $Zg$.
\end{lemma}

\begin{proof}
Write $g = \qq(A,B,\lambda)$.
We would like to solve, for $T \in N$ and
$S \in Z$, the equation
$$
\nn(T) \qq(A,B,\lambda) \nn(-T) = S\qq(A,B,\lambda),
$$
or,
$$
\left(\begin{matrix}I &T \\0 &I \end{matrix}\right)
\left(\begin{matrix}A  & 0  \\0 & \lambda A^* \end{matrix}\right)
\left(\begin{matrix}I & B \\0 &I \end{matrix}\right)
\left(\begin{matrix}I & -T \\0 &I \end{matrix}\right) =
\left(\begin{matrix}I & S \\0 &I \end{matrix}\right)
\left(\begin{matrix}A  & 0 \\0 & \lambda A^*  \end{matrix}\right)
\left(\begin{matrix}I & B \\0 &I \end{matrix}\right),
$$
or,
$$
\left(\begin{matrix}I &T \\0 &I \end{matrix}\right)
\left(\begin{matrix}A  & 0  \\0 & \lambda A^* \end{matrix}\right)
\left(\begin{matrix}I & -T \\0 &I \end{matrix}\right) =
\left(\begin{matrix}I & S \\0 &I \end{matrix}\right)
\left(\begin{matrix}A  & 0 \\0 & \lambda A^*  \end{matrix}\right),
$$
or,
$$
\lambda T A^* - AT = \lambda S A^*.
$$
Observe that for matrices $L_1$ and $L_2$, the transformation
$$
X \mapsto L_1 X - X L_2
$$
is singular if and only if an eigenvalue of $L_1$ is the same
as an eigenvalue of $L_2$. This implies that the equation
$$
\lambda T A^* - AT = \lambda S A^*
$$
can be solved for $T$ if $A$ and $\lambda A^*$ do not share an eigenvalue,
i.e., if the eigenvalues of $A$ are $\{\beta_1, \beta_2 \}$,
then
$$
\{\beta_1, \beta_2 \} \cap \lambda\{\beta^{-1}_1, \beta^{-1}_2 \}
= \phi,
$$
i.e., $\lambda \not \in \{\beta_1^2, \beta_1 \beta_2, \beta_2^2 \}.$

We actually need to solve for $T$ with
$$
\omega {}^t\!T \omega = T.
$$
For this we write the earlier equation as
$$
\lambda T  - ATA^{*-1} = \lambda S.
$$
If $ T \mapsto \omega {}^t\!T \omega $ is denoted by $\sigma$,
then the above equation becomes:
$$
\lambda T  - AT\sigma(A) = \lambda S.
$$
Applying $\sigma$ to this equation, we obtain
$$
\lambda \sigma(T)  - A \sigma(T) \sigma(A) = \lambda S.
$$
Adding the two previous equations,
$$
\lambda [T+ \sigma(T)]  - A[T + \sigma(T)]\sigma(A) = 2\lambda S.
$$
Now $\frac12(T + \sigma(T))$ is of the desired form, completing the
proof of the lemma.
\end{proof}

Suppose from now on that the hypothesis
of lemma~\ref{lem:eigens} does not hold.
Then we can divide the rest of the proof into three cases.

\textbf{Case 1}:
$\beta_i^2=\lambda$ for precisely one value of $i\in\{1,2\}$.
Assume without loss of generality that
$\beta_1^2 = \lambda \neq \beta_2^2$.
Since $\beta_1$ and $\beta_2$ thus have different minimal polynomials
over $k$, they must both lie in $k\mult$.
Therefore, for $A' = A/\beta_1$, $A'^*=\beta_1A^*
= (\lambda A^*)/\beta_1$. Thus for $g' = g/\beta_1$, the $2 \times 2$
block diagonal matrices are $(A',A'^*)$, i.e., up to scaling $g$
belongs to $\Sp(4)$. Furthermore, one of the
eigenvalues of $A'$ is 1. By appealing to \cite{br}, we will
see in \S\ref{sec:completed} that we don't have to worry about
these elements.

\textbf{Case 2}:
$\lambda = \beta_1^2 = \beta_2^2$.
Then $\beta_2 = \pm \beta_1$.
If $\beta_2=\beta_1$, then
$\lambda = \beta_1\beta_2 = \det(A)$
and $\tr(A) \neq 0$.
If $\beta_2 = -\beta_1$,
then $\lambda = -\beta_1\beta_2 = -\det(A)$
and $\tr(A) = 0$.

\begin{lemma}
\label{lem:tr-notzero}
If $\tr(A) \neq 0$, then the $J$-orbit of $g= \qq(A,B,\det A)$
is $\tau$-invariant.
\end{lemma}

\begin{proof}
In this case, $u=-\lambda$.
Let $n = \nn\left(\begin{smallmatrix}x&y\\z&x\end{smallmatrix}\right)$,
where $x,y,z\in k$ are to be determined.
Then
$$
n g n\inv =
\qq(
\left(\begin{smallmatrix}au&1 \\u&0\end{smallmatrix}\right),
\left(\begin{smallmatrix}R&S \\T&R\end{smallmatrix}\right),
\lambda),
$$
where $(R,S,T)$ is an affine function of $(x,y,z)$
that takes the value $(r,s,t)$ at the origin
and has gradient
$$
\begin{pmatrix}
\lambda u\inv -1 & 0 & -a\lambda u\inv \\
0 & \phnm -1 & \lambda u^{-2} \\
-2a\lambda & \lambda & a^2\lambda - 1
\end{pmatrix}.
$$
Since $a\neq 0$, the first and third rows are independent
(look at the second and third columns).
Thus, we may choose $x$, $y$, and $z$ to give
$R$ and $T$ any desired value, so $g$ is $N$-conjugate
(and thus $J$-conjugate) to
$$
g' =
\qq(
\left(\begin{smallmatrix}-a\lambda &1 \\-\lambda &0\end{smallmatrix}\right),
\left(\begin{smallmatrix}0&s' \\0&0\end{smallmatrix}\right),
\lambda)
$$
for some $s'\in k$.
Let
$$
p =
\hh(0,-s')
\:
\hh((a,0),0)
\:
\mm(1,\left(\begin{smallmatrix}\phm0&1\\-1&0\end{smallmatrix}\right))
\in J.
$$
Then
$pg'p\inv = \tau(g')$.
Thus, the $J$-orbit of $g$ is $\tau$-invariant.
\end{proof}

\begin{lemma}
\label{lem:tr-zero}
If $\tr(A)=0$, then the $J$-orbit of $g =\qq(A,B,\pm\det A)$
is $\tau$-invariant.
\end{lemma}

\begin{proof}
Let $h=\hh(0,-t\lambda\inv)$.  Let
$$
m =
\begin{cases}
\mm(1,\left(\begin{smallmatrix}0&-1\\1&\phm0\end{smallmatrix}\right))
& \text{if $\lambda = -\det(A)$}, \\
\mm(1,\left(\begin{smallmatrix}\phm 0&1\\-1&0\end{smallmatrix}\right))
& \text{if $\lambda = \phm\det(A)$}.
\end{cases}
$$
Let $g'=hgh\inv$.
Then
$$
g' = \qq(
	\left(\begin{smallmatrix}0&1\\u&0\end{smallmatrix}\right),
	\left(\begin{smallmatrix}r&s'\\0&r\end{smallmatrix}\right),
	\lambda)
$$
for some $s'\in k$.
Let $h'=\hh(0,-s')$.
Then
$h'm\tau(g')m\inv {h'}\inv = g'$.
Thus, the $J$-orbit of $g$ is $\tau$-invariant.
\end{proof}

\textbf{Case 3}:
$\lambda \neq \beta_i^2$ for $i=1,2$.
Then $\lambda = \beta_1\beta_2 = \det(A)$
(since we are assuming that the hypothesis of
lemma~\ref{lem:eigens} is not satisfied).
Therefore either lemma~\ref{lem:tr-notzero}
or lemma~\ref{lem:tr-zero} applies.

\section{Closed cell}
\label{sec:closed}

In this section we prove the following:
\begin{prop}
\label{prop:closed}
Any distribution on $P$ that is
$J$-invariant and $\psi$ invariant for a non-trivial
character $\psi$ of $Z$
must be invariant under the involution $\tau$.
\end{prop}
\begin{proof}
Let $p = \mm(\lambda, A)h$ be an element of $P=MH$.
It is easy to see that for $z \in Z$,
$$
pzp^{-1}= (\lambda^2/\det A) z.
$$
Therefore if  $\lambda^2 \neq \det A$, then for any $z_0 \in Z$, there is
$z \in Z$ such that
$$
pzp^{-1} = zz_0,
$$
or, $z^{-1}pz = z_0p,$ implying that the $J$ orbit of such a $p$
is $Z$ stable.

On the other hand, if $\lambda^2 = \det A$, then the
$J$-orbit of $p$ is $\tau$-invariant.
This can be either
checked as an easy exercise, or else observe that in this case
$p/\lambda$ in fact belongs to $\Sp(4)$, and therefore one can use
the calculation of lemma~5.4 of \cite{br}.
\end{proof}

\section{Constructible sets}

Before we put all of the pieces together to prove our main theorem,
we need a bit of general
topology that does not seem to have
been carefully written down anywhere that we could find.
The reason for our need is that
a distribution can be restricted to an open set,
and we can try to decide
if the distribution is zero or not on it.
If zero, then the question becomes one on
the complementary closed set. And we can proceed inductively trying
to prove that a distribution is zero on the whole set, the kind of goal
we have set ourselves to in this paper.

However, situations might arise where a space is decomposed not into
an open and a complementary closed set, but into a slightly more complicated
subset, a \emph{constructible} set, and its complement.
In our case this arises when we are considering all $J$-conjugates
of elements for which the hypothesis of lemma~\ref{lem:eigens}
holds, where we would like to apply lemma~\ref{lem:Z-stable}
to say that the distribution restricted to such
elements is zero, except that it does not make sense to restrict
distributions to such general subsets (a certain union of orbits).

First let us recall that a subset $Y$ of a topological space $X$ is said
to be \emph{constructible} if $Y$ is a finite union of locally closed subsets.
(A subset is \emph{locally closed}
if it is the intersection of a closed set
with an open set.)

The reason for the importance of constructible sets in $p$-adic groups
arises from the following theorem, which is a variation of a
well-known theorem due to Chevalley in algebraic geometry. We refer
to \cite{bz} for a proof.

\begin{thm}
\label{thm:bz-constructible}
Let $\underline{X}$ and $\underline{Y}$ be algebraic
varieties over a non-archimedean local field $k$,
and $\underline{f}$ be a morphism
of algebraic varieties between $\underline{X}$ and $\underline{Y}$.
Denote the corresponding $k$-valued points, and the morphism between
the $k$-valued points by removing the underline. Then $f(X)$ is a
constructible subset of the topological space $Y$.
\end{thm}

For our purposes, the following lemma is of utmost importance.

\begin{lemma}
\label{lem:constructible}
Let $Y$ be a constructible subset of a topological space
$X$. Then there are finitely many closed subsets
$X_1\subset \cdots \subset X_n = X$
such that $X_{i+1}\smallsetminus X_i$
is an open subset of $X_{i+1}$ which is either
contained in $Y$ or in $X \smallsetminus Y$. Further, this decomposition
is \emph{canonical} in the sense that if a group
operates on $X$ preserving $Y$, then it also preserves
each of the sets $X_i$.
\end{lemma}

\begin{proof}
We recall the following well-known decomposition of
a constructible set, cf.\ \cite{bz},
into a disjoint union of locally closed subsets.

For any subset $A$ of $X$, define
$C^1(A) = C(A) = \bar{A} \smallsetminus \overline{ \bar{A}\smallsetminus A}$.
Clearly, $C(A)$ is a locally closed subset of $X$, and is contained in $A$. 
Further, $C(A)=A$ if and only if $A$ is locally closed.
For $i>1$, inductively define $C^i(A)$ to be
$C(A\smallsetminus [C(A)\cup C^2(A) \cup \cdots \cup C^{i-1}(A)])$.
It is easy to see
that if $A$ is constructible, then $C^i(A)$ is empty for large $i$,
and therefore such an $A$ is a finite \emph{disjoint} union of
the locally closed sets $C^i(A)$.

Renaming the indices, let $Y = \cup_{i=1}^{n-2} Y_i$,
a disjoint union of locally closed subsets
$Y_i$ with $Y_i = Z_i \smallsetminus W_i$ where $Z_i$ and $W_i$ are
closed subsets of $X$.
Now define $X_1 = \cup_{i=1}^{n-2}W_i$,
$X_i = X_{i-1} \cup Z_{i-1}$ for $1<i<n$,
and $X_n = X$.

Clearly the $X_i$'s have the property desired. These sets are
canonically constructed, and therefore are preserved
under any group action preserving $Y$.
\end{proof}

\section{Proof of uniqueness of Fourier-Jacobi models completed}
\label{sec:completed}

We now have all the pieces necessary
to complete the proof of theorem~\ref{thm:T-zero}.
We start with a distribution $T$ with invariance properties as in the
statement of this theorem.
Our aim is to prove
that such a distribution is identically zero. By appealing to
\cite{br} as in \S\ref{sec:br}, we already know that $T$ is zero
on ${\mathcal G} = k\mult\Sp(4)$.
By \S\ref{sec:open}, $T$ is zero on the open cell,
thus $T$ is supported on the union $P \cup Pw'P$ of the closed and
the middle cell.
Let $Y=Pw'P$, an open subset of this union.
Write $Y = \Yo \cup \Yc$, with $\Yo$ the (open) subset
of $Y$ on which $\blah(g) \not = 0$.

By lemmas~\ref{lem:middle-open} and \ref{lem:bernstein},
$T$ is zero on $\Yo$, thus $T$ is supported
on $P\cup \Yc$. Since anyway we know that the support of $T$ is outside
${\mathcal G}$,
we get that the support of $T$ is contained in the closed subset
$\Ycc= \Yc \smallsetminus (\Yc \cap {\mathcal G})$ of $\Yc$.
Let $S$ denote the set of elements of the form appearing
in lemma~\ref{lem:eigens}.
This is the set of rational points of a $k$-variety.
One can write
$\Ycc = Y_1 \cup Y_2$, where $Y_1$ is the subset
of $\Ycc$ for which the hypothesis of lemma~\ref{lem:eigens} holds,
i.e., it consists of $J$-conjugates of elements of $S$.
Applying theorem~\ref{thm:bz-constructible}
to the map from $S \times J$ to $\GSp(4)$
defined by $(s,j)\mapsto jsj^{-1}$,
we see that
$Y_1$ is a constructible
subset of $\Ycc$.
Applying lemma~\ref{lem:constructible}
to the toplogical space $X = \Ycc$,
and $Y = Y_1$, we are able to write $\Ycc$ as an increasing
union of closed sets such that the successive differences are
either in $Y_1$ or $Y_2$, to which we can apply now
lemmas~\ref{lem:Z-stable} and~\ref{lem:bernstein}
respectively (for the first, $J$ orbits are $Z$-invariant, and for the
second, $J$-orbits are $\tau$-invariant by
lemmas~\ref{lem:tr-notzero} and~\ref{lem:tr-zero})
to conclude that the distribution is zero on $\Ycc$,
thus is supported in $P$.

We remind the reader that by removing
$\Yc \cap {\mathcal G}$ from $\Yc$, we have removed from our
consideration
elements with eigenvalues
$\{c ,c \alpha, c\alpha\inv,c \}$,
where $\alpha,c\in k\mult$ and $\alpha \neq \pm 1$, for which $J$-orbits are
in fact not $\tau$-invariant, and were a source of difficulty
for \cite{br}; for us, luckily, we can just use \cite{br}
instead of having to redo this part of their argument.

Finally, by proposition~\ref{prop:closed}, $T$ is zero on $P$,
completing the proof of theorem~\ref{thm:T-zero}.

\section{Proof of multiplicity one theorem for $\GSp(4)$ for
nonsupercuspidals}
\label{sec:non-sc}

\begin{prop}
Any irreducible admissible nonsupercuspidal representation
$\tpi$ of $\GSp(V)$ decomposes with multiplicity one when
restricted to $\Sp(V)$;
i.e., for any irreducible, admissible representation
$\pi$ of $\Sp(V)$,
$$
\dim \Hom_{\Sp(V)}[\tpi, \pi] \leq 1.
$$
\end{prop}

\begin{proof}
From generalities, it is known that $\tpi$ decomposes as a finite
direct sum of irreducible representations when restricted to $\Sp(4)$,
say
$$
\tpi|_{\Sp(4)} \cong \pi_1\oplus \cdots \oplus \pi_\ell,
$$
where we allow repetitions. One has the relation \cite{gk}:
$$
\dim \End_{\Sp(4)}(\tpi|_{\Sp(4)})
=
\card\set{\chi:k\mult \rightarrow \C\mult}{ \tpi \chi \cong \tpi},
$$
where $\chi$ is a character of $k\mult$ thought of as a character of
$\GSp(2n)$, via the similitude homomorphism
$\lambda: \GSp(2n) \rightarrow k\mult$. For a character $\chi$ of $k\mult$,
and a representation $\tpi$ of $\GSp(4)$, the representation
$\tpi \chi$ is said to be a twist of $\tpi$.
The characters $\chi$ of $k\mult$
such that $ \tpi \chi \cong \tpi$ are called self-twists
of $\tpi$. Observe that the above relation implies that if the
number of self-twists of $\tpi$ is less than $4$,
then $\tpi$
decomposes with multiplicity one when restricted to $\Sp(4)$.

Any nonsupercuspidal
representation $\tpi$ of $\GSp(2n)$ can be realized
as a subquotient of a
representation induced from a supercuspidal representation of a
Levi factor of a parabolic subgroup of $\GSp(2n)$.
Recall that this \emph{cuspidal support}
is uniquely determined by the representation $\tpi$ (up to
conjugation by the normalizer of the Levi subgroup). In particular,
if $\tpi$ has self-twists, then so does the associated
cuspidal support (where twisting is defined in the obvious way),
up to equivalence.

Recall that the parabolic subgroups in $\GSp(2n)$ are parametrized
by partitions $n=n_1+\cdots + n_r + m$ with corresponding Levi
subgroup isomorphic to $\GL(n_1) \times \cdots \times \GL(n_r)
\times \GSp(2m)$. (In this notation, $m=0$ is allowed; $\GSp(0)$ is taken
to be $k\mult$.) Irreducible representations of such Levi subgroups are
thus parametrized by $(\tau_1,\cdots, \tau_r; \tau)$ where $\tau_i$
are irreducible admissible representations of $\GL(n_i)$, and $\tau$
is one of $\GSp(2m)$. Further,
the twisting operation
is given by
$$
(\tau_1,\cdots, \tau_r; \tau)\chi = (\tau_1,\cdots, \tau_r; \tau\chi).
$$

We take up the various possibilities for the cuspidal support
of $\tpi$.

\textbf{Case 1 (Siegel Parabolic):}
In the notation above, this corresponds to $n_1=2$,
and $m=0$.
The Levi subgroup $M$ is isomorphic to $\GL(2) \times k\mult$.
Denote a representation of $M$ by $(\tau;\chi)$ with $\tau$ a supercuspidal
representation of $\GL(2)$, and $\chi$ a character of $k\mult$.
The action of the
normalizer of $M$ on this data is given by
$$
(\tau; \chi) \mapsto (\tau^\vee; \omega_\tau\chi),
$$
where $\omega_\tau$ is the central character of $\tau$.
Since the twisting is given by $(\tau;\chi) \nu= (\tau;\chi \nu)$,
it follows that the cuspidal support has a non-trivial self-twist
if and only if
$$
(\tau; \chi \nu) = (\tau^\vee;  \omega_\tau\chi),
$$
i.e., if and only if $\tau \cong \tau^\vee$, and $\nu = \omega_\tau$.
Clearly there is at most 1 non-trivial self-twist, therefore
$\tpi$ splits into at most two irreducible factors when restricted to
$\Sp(4)$. These factors must be distinct.

\textbf{Case 2 (Klingen Parabolic):}
In the notation introduced above,
this corresponds to $n_1 = m =1$.
The Levi subgroup $M$ is isomorphic to $k\mult \times \GL(2)$.
Denote a representation of $M$ by $(\chi; \tau)$ with $\tau$ a supercuspidal
representation of $\GL(2)$, and $\chi$ a character of $k\mult$.
The action of the
normalizer of $M$ on this data is given by
$$
( \chi; \tau ) \mapsto (\chi^{-1}; \chi \tau).
$$
Since the twisting is given by $(\chi; \tau) \nu= (\chi; \tau\nu)$,
it follows that the cuspidal data has a self-twist by $\nu$
if and only if either of the following holds:
\begin{align*}
(\chi;\tau) & \cong  ( \chi;\tau\nu) ,\\
(\chi^{-1};\chi \tau) & \cong  (\chi; \nu \tau).
\end{align*}
Thus $\nu$ is a self-twist for $\tpi$
if and only if either $\nu$ is a
self-twist for $\tau$, or  $\chi^2=1$, and $\chi \nu$
is a self-twist for $\tau$. Thus if ${\mathcal S}_\tau$ denotes the
set of self-twists of $\tau$, then the self-twists of the cuspidal pair
$(\chi; \tau)$ are parametrized by elements of order 2 in
$ {\mathcal S}_\tau \cup \chi
{\mathcal S}_\tau$, and thus the
cardinality of the set of self-twists of the cuspidal pair
$(\chi; \tau)$ is at most $2s$, where $s$ is the cardinality
of ${\mathcal S}_\tau$. By multiplicity one restriction
from $\GL(2)$ to $\SL(2)$, we already know that there are
$s$ irreducible components in the restriction
of $\tau$ to $\SL(2)$. Let $\pi_1$ be an irreducible representation
of $\Sp(4)$ appearing inside $\tpi$.
Clearly, the cuspidal support
of $\pi_1$ (a representation of $\Sp(4)$) is $(\chi; \tau_1)$ where
$\tau_1$ is an irreducible representation of $\SL(2)$ appearing in the
restriction of $\tau$.
Since all the
conjugates of $\pi_1$ by $\GSp(4)$ also appear in the restriction
of $\tpi$ to $\Sp(4)$, we find irreducible representations
of $\Sp(4)$ inside $\tpi$
with cuspidal support $(\chi;\tau_i)$ where $\tau_i$ are the various
irreducible representations of $\SL(2)$ inside $\tau$.
This gives us $s$ irreducible representations
of $\Sp(4)$ inside $\tpi$.

It is clear that each irreducible
representation of $\Sp(4)$ appears in $\tpi$ with the same multiplicity.
If this common multiplicity were more than $1$,
then each irreducible component would contribute at least $4$
to the dimension of the endomorphism ring of $\tpi$ restricted to $\Sp(4)$.
But this would contradict the fact that
the dimension of the endomorphism ring
is bounded by $2s$.

\textbf{Case 3 (The Borel):}
In the notation above, this corresponds to $n_1=n_2=1$,
$m=0$. The Levi subgroup $M$ is isomorphic to $k\mult \times k\mult \times
k\mult$.
Denote a representation of $M$ by $(\chi_1,\chi_2; \chi)$
where $\chi_i$ and $\chi$ are characters of $k\mult$.
The relevant Weyl group in this case is the semi-direct product
of $(\Z/2)^2$ by $\Z/2$; its action on the
character $(\chi_1,\chi_2; \chi)$ of  $M$ is given by
compositions of the actions:
\begin{align*}
(\chi_1,\chi_2; \chi) &\mapsto (\chi_2,  \chi_1; \chi), \\
(\chi_1,\chi_2; \chi) &\mapsto (\chi_1^{-1},  \chi_2\phinv ; \chi \chi_1 ), \\
(\chi_1,\chi_2; \chi) &\mapsto (\chi_1\phinv,  \chi_2^{-1} ; \chi \chi_2 ).
\end{align*}
In this case, twisting is given by
$(\chi_1,\chi_2; \chi)\nu = (\chi_1,\chi_2; \chi\nu)$.
From these facts, it can be checked that
the self-twists on the cuspidal data
$(\chi_1,\chi_2; \chi)$
are the
characters belonging to the group generated by characters of order
2 in $\{\chi_1,\chi_2 \}$, a group of order $1$, $2$, or $4$.
If the number of self-twists is 1 or 2, there is nothing further
to be said. If the number of self-twists is 4, observe that
$\chi_1$ and $\chi_2$ must be both of order 2, and therefore
by lemma~3.2 of \cite{st}, the principal series representation
induced from the cuspidal data is irreducible
as a representation of $\GSp(4)$, and we are done by the uniqueness
of Whittaker models.
\end{proof}

\section{Proof of multiplicity one theorem for $\GSp(4)$ completed:
existence of Fourier-Jacobi models for supercuspidal representations}
\label{sec:existence}

The multiplicity one theorem for $\GSp(4)$ will be completed
(by this method of Fourier-Jacobi models;
we already have given an independent proof in \S\ref{sec:mult1})
if we can show that every
representation of $\GSp(4)$ has a Fourier-Jacobi model.
Unlike the Whittaker model, which may not
exist for some representations, Fourier-Jacobi models (should)
always exist
as long as $\pi$ is not one dimensional.
Unfortunately we are able to prove this
only for supercuspidal representations of $\GSp(4)$;
this is enough for our purposes,
as we have already proved
multiplicity one in \S\ref{sec:non-sc}
for all nonsupercuspidal representations.

We begin with the following general lemma.

\begin{lemma}
Let $\pi$ be a smooth representation of
$N\cong k$, on which $N$ acts nontrivially.
Then there exists a
nontrivial additive character $\psi: N \rightarrow \C\mult$,
and a nonzero linear form $\ell: \pi \rightarrow \C$
such that
$$
\ell(nv) = \psi(n) \ell(v),
$$
for all $v \in \pi, n \in N$.
\end{lemma}

\begin{proof}
It clearly suffices to prove that the twisted Jacquet module
$$
\pi_\psi := \frac{\pi}{ \set{ nv-\psi(n)v}{n \in N, v \in \pi}}
$$
is nonzero for some
nontrivial additive character $\psi: N \rightarrow \C\mult$.
This will be a simple consequence of the exactness of the Jacquet
functor, denoted $\pi \mapsto \pi_N$
defined as above but for $\psi =1$, and of the twisted Jacquet
functor $\pi \mapsto \pi_\psi$. Let $\pi[N]$ be the kernel of
the map from $\pi$ to $\pi_N$. Clearly $\pi[N]$ is an $N$-module,
which by the exactness of the Jacquet functor has trivial
Jacquet module. Since any finitely-generated representation has an
irreducible quotient, cf.\ \cite{bz} lemma~2.6(a), any smooth
representation has an irreducible \emph{subquotient}. By Schur's lemma,
which is valid
for any smooth representation, cf.\ \cite{bz} lemma~2.11, any irreducible
representation of $N$ is one dimensional.
Therefore $\pi[N]$ has an irreducible
subquotient which is one dimensional,
and therefore given by a non-trivial character
$\psi: N \rightarrow \C\mult$.
The exactness of the twisted Jacquet functor then implies
that the twisted
Jacquet module of $\pi[N]$, and hence of $\pi$, is non-trivial.
\end{proof}

\begin{lemma}
Let $\pi$ be an irreducible
smooth representation of $\GSp(4)$
that is not one dimensional.
Let $\psi$ be a nontrivial character of $k$, thought of as a
character of $Z$, the center of the unipotent radical of the
Klingen parabolic subgroup of $\GSp(4)$. Then
$$
\Hom_{Z}(\pi, \psi) \not = 0.
$$

\end{lemma}

\begin{proof}
It is easy to see that a Levi subgroup of the Klingen parabolic operates
transitively on the set of all non-trivial characters of $Z$.
Hence,
if the conclusion of the lemma is true for one $\psi$, it is true for
all $\psi$.
The previous lemma gives that if
$\Hom_{Z}(\pi, \psi) = 0$,
then $Z$ acts trivially on $\pi$, and hence
the normal subgroup generated by $Z$ also acts trivially on $\pi$.
But it is a standard fact that $\Sp(4)$ has no normal subgroup besides
the center, and hence $\Sp(4)$ must act trivially on $\pi$, and therefore
$\pi$ must be one dimensional, concluding the proof of the lemma.
\end{proof}

We still have the task of proving that the representation
$$
\pi_\psi := \frac{\pi}{ \set{ zv-\psi(z)v}{z \in Z, v \in \pi}},
$$
of $J$, which we know now is nonzero, has nonzero
irreducible quotients.
We would have liked to believe that this is obvious, but we did
not succeed in finding a general proof.
Here is a proof for the case where $\pi$ is supercuspidal.

In the next two lemmas, $G$ is a general $\ell$-group, countable at
infinity, in the sense of \cite{bz}.
This hypothesis is satisfied
by algebraic groups over non-archimedean local fields.
We let ${\mathcal S}(G)$ denote the Schwartz space of locally constant,
compactly supported functions on $G$, thought of as a left $G$-module.
We let $dg$ denote a Haar measure on $G$.

We recall proposition~2.12 of \cite{bz}:

\begin{lemma}
Let $G$ be an $\ell$-group, and $f$ a compactly supported
function on $G$. Then there is an irreducible smooth representation
$\pi$ of $G$ such that the action of $f$ on $\pi$ is non-trivial.
\end{lemma}

We combine this lemma with the following trivial lemma:

\begin{lemma}
Let $\pi$ be a smooth irreducible representation of $G$. Then for
every vector $v \in \pi$, there is a homomorphism of $G$-modules
${\mathcal S}(G) \rightarrow  \pi$ given by
$$
f \mapsto \int_G f(g) \pi(g) v \, dg.
$$
For a function $f \in {\mathcal S}(G)$,
the image of $f$ under this homomorphism is non-zero for
some choice of $v \in \pi$ if and only if the action of $f$
on $\pi$ is nontrivial.
\end{lemma}

\begin{prop}
Let $\pi$ be a supercuspidal representation of $\Sp(4)$.
Then the representation $\pi_\psi$ of $J$ has an irreducible quotient.
\end{prop}

\begin{proof}
Observe that the supercuspidal representation
$\pi$ can be realized on a space of functions in
${\mathcal S}(\Sp(4))$. Fix one such realization, and think of
elements of $\pi$ now as functions on $\Sp(4)$.
Restricting
these functions to the Fourier-Jacobi group $J$, we get a
space of locally constant, compactly supported functions on $J$.
For a function $g$ of this kind, and for any
element $z \in Z$,  $f = zg-g$ is another
such function.
We can (and do) choose $z\in Z$ so that $f$ is nonzero.
By the previous two lemmas,
there is an irreducible representation $\rho$ of $J$
on which $f$ acts nontrivially.
By generalities (cf.\ \cite{bz}, proposition~2.11),
$\rho$ has a central character
(i.e., Schur's lemma holds).
Hence, $Z$ operates by a character on $\rho$.
This character cannot be
trivial, as $f$ was chosen to be of the form $zg-g$.
\end{proof}

\begin{question}
It would be interesting to understand $\pi_\psi$
as a representation of $J$.
Of the
irreducible representations of $J$ with central character $\psi$
(which are parametrized by irreducible representations of $\tilde{S}$,
the two-fold cover of $\SL(2)$),
which ones occur as a quotient in $\pi_\psi$?
We expect that if $\pi$ is a generic representation of
$\GSp(4)$, then every irreducible representation of $\tilde{S}$
appears.
Further, if $\pi$ is a degenerate representation, then
we expect $\pi_\psi$ to be a representation of finite length as a
$J$-module, with a unique irreducible quotient.
\end{question}

\end{document}